\documentclass[12pt,letterpaper]{article}
\usepackage{fullpage}

\usepackage{amsmath,amsthm,amssymb,amsfonts}
\renewcommand{\Bbb}{\mathbb}
\usepackage{mathtools}
\usepackage{mathrsfs}
\usepackage[shortlabels]{enumitem}

\usepackage{comment}
\usepackage{hyperref}
\usepackage[nameinlink,capitalise,noabbrev]{cleveref}
\crefformat{equation}{(#2#1#3)}

\usepackage{xcolor}
\usepackage{graphicx}
\usepackage{float}
\usepackage{tikz}
\usepackage{pgfplots}
\pgfplotsset{compat=1.18}

\theoremstyle{definition}
\newtheorem{definition}{Definition}

\theoremstyle{plain}
\newtheorem{conjecture}[definition]{Conjecture}
\newtheorem{theorem}[definition]{Theorem}
\newtheorem{lemma}[definition]{Lemma}

\newtheorem{claim}[definition]{Claim}

\crefname{claim}{Claim}{Claims}
\crefname{lemma}{Lemma}{Lemmas}

\def \l {\ell}

\def \sm {\setminus}
\def \ce {\coloneqq}
\def \E {\mathbb{E}}
\def \P {\mathbb{P}}
\def \F {\mathbb{F}}

\renewcommand{\le}{\leqslant}
\renewcommand{\ge}{\geqslant}
\renewcommand{\leq}{\leqslant}

\makeatletter
\def \eps {\varepsilon}
\def \es {\varnothing}
\renewcommand \b[2] {\binom{#1}{#2}}

\def \A {\mathcal{A}}
\def \B {\mathcal{B}}

\def \C {\mathcal{C}}
\def \mF{\mathcal{F}}
\def \G {\mathcal{G}}
\def \mH{\mathcal{H}}
\def \mP{\mathcal{P}}

\def \S {\mathcal{S}}
\def \K {\mathcal{K}}

\newcommand{\ka}[1][A]{K_{#1}}
\newcommand{\kca}[1][A]{K^\circ_{#1}}
\newcommand{\nka}[1][A]{n_{K_{#1}}}
\newcommand{\nkca}[1][A]{n_{K^\circ_{#1}}}

\newcommand{\hu}[1][U]{\mH_{#1}}
\newcommand{\ha}[1][A]{\mH_{#1}}
\newcommand{\fa}[1][A]{\mF_{#1}}

\newcommand{\eha}{\E[|\ha|]}
\newcommand{\ehu}{\E[|\hu|]}

\newcommand{\dfas}[1][S]{d_{\mF_{A}}(#1)}
\newcommand{\dhas}[1][S]{d_{\mH_{A}}(#1)}
\newcommand{\dhus}[1][S]{d_{\mH_{U}}(#1)}

\newcommand{\edhas}[1][S]{\E[d_{\mH_{A}}(#1)]}
\newcommand{\edhus}[1][S]{\E[d_{\mH_{U}}(#1)]}

\newcommand{\Onoverq}{O \left( \frac{n}{q}\right) }
\newcommand{\Omenoverq} {\Omega \left( \frac{n}{q}\right) }

\newcommand{\Omenpoverq}[1]{\Omega \left( \frac{n^{#1}}{q}\right) }

\newcommand{\alFqd}[1][d]{\alpha(\F_q^{#1})}
\newcommand{\alFqdp}[1][d]{\alpha(\F_q^{#1},p)}

\newcommand{\alFthrdp}{\alpha(\F_q^{3},p)}

\title{Maximum number of points in general position in a random subset of finite $3$-dimensional spaces}

\author{    J\'ozsef Balogh\footnote{
    Department of Mathematics, University of Illinois at Urbana--Champaign, Urbana, IL, USA. Research supported in part by NSF grants RTG DMS-1937241, FRG DMS-2152488, the Arnold O.~Beckman Research Award (UIUC Campus Research Board RB 24012), and the Simons Fellowship. E-mail: \texttt{jobal@illinois.edu}.} \and
    Haoran Luo\footnote{
    Department of Mathematics, University of Illinois at Urbana--Champaign, Urbana, IL, USA. Research supported in part by NSF grants RTG DMS-1937241, FRG DMS-2152488, the Arnold O.~Beckman Research Award (UIUC Campus Research Board RB 24012) and Dr. James J. Woeppel Fellowship, E-mail: \texttt{haoranl8@illinois.edu.}}
}
\date{}

\begin{document}
\maketitle

\begin{abstract}
    Let $\alpha(\mathbb{F}_q^{d},p)$ be the maximum possible size of a point set in general position in the $p$-random subset of $\mathbb{F}_q^d$. In this note, we determine the order of magnitude of $\alpha(\mathbb{F}_q^{3},p)$ up to a polylogarithmic factor by proving a balanced supersaturation result for the sets of $4$ points in the same plane.
\end{abstract}

\section{Introduction} \label{sec::int}
For a field $\F$ and a positive integer $d$, let $\F^d$ be the $d$-dimensional affine space over $\F$.
A point set $S \subseteq \F^d$ is \emph{in general position} if there are no $d+1$ points in $S$ contained in a $(d-1)$-dimensional affine subspace (a \emph{$(d-1)$-flat}). This notion is closely related to many central problems in Discrete Geometry~\cite{dudeney1959amusements, erdos1986some, furedi1991maximal, cardinal2017general, milicevic2017sets, balogh2018number, suk2023higher, nenadov2024number}, and we refer the reader to Chapter 8 in~\cite{brass2005research} for a survey.

In this note, we focus on the case where $\F = \F_q$, the finite field with $q$ elements. Let $\alFqd$ be the maximum possible size of a point set in general position in $\F_q^d$. Erd\H{o}s noticed that $\alFqd \ge q$ by the construction $\{(x,x^2,\ldots, x^d): x \in \F_q\}$, see Appendix in~\cite{roth1951problem}. For the upper bound, it is easy to see that $\alFqd \le dq$, as  $\F_q^d$ can be partitioned into $q$ disjoint $(d-1)$-flats, and each can contain at most $d$ points. A more careful argument shows that $\alFqd \le (1+o(1))q$,
see for example Lemma~4.1 in~\cite{chen2023random}.
Roche-Newton and Warren~\cite{rochenewton2022arcs} started to work on  a random version of this problem. Let $\alFqdp$ be the maximum possible size of a point set in general position in a $p$-random subset of $\F_q^d$,  where \emph{$p$-random} means that every point is chosen with probability $p$ independently. They~\cite{rochenewton2022arcs} gave essentially tight bounds on $\alFqdp[2]$ for $p\in [0,q^{-1}] \cup [q^{-1/3},1]$. Recently, Chen, Liu, Nie, and Zeng~\cite{chen2023random} determined the order of magnitude of $\alFqdp[2]$ up to a polylogarithmic factor for all possible values of $p$ and the order of magnitude of $\alFqdp$ when $p$ is in certain ranges.

For the general case, the following conjecture is natural, based on the previous results mentioned above. It was also implicitly mentioned in Section~6 in~\cite{chen2023random}.

\begin{conjecture} \label{conj::generalTheorem}
    For every integer $d \ge 2$, there is a  positive real number $C = C(d)$ such that the following is true. As the prime power $q$ goes to infinity, we have with high probability
    \[
        \alFqdp =
        \left\{
        \begin{array}{ll}
         \Theta(pq^d)     &  \textrm{ for } \quad \omega(q^{-d}) = p = o(q^{-d + 1/d}), \\
         \Theta(pq)    & \textrm{ for  } \quad Cq^{-1+1/d} \log^2 q \le p  \le 1,
        \end{array}
        \right.
    \]
    and
    \[
        \Omega(q^{1/d}/\log q) = \alFqdp = O(q^{1/d}\log^2q) \quad \textrm{ for }\quad  q^{-d+1/d}/\log q \le p < Cq^{-1+1/d}\log^2q,
    \]
    where all the implicit constants depend only on $d$.
\end{conjecture}

We will show that \cref{conj::generalTheorem} follows from the following conjecture about balanced supersaturation.
For every set $U \subseteq \F_q^d$ and every integer $i \in \{1,2,\ldots, d+1\}$, let $\G_{i}(U)$ be the collection of $i$-subsets of $U$ that are not in any $(i-2)$-flat and $\S_{i}(U) \ce \b{U}{i} \sm \G_{i}(U)$.
For every $r$-uniform hypergraph $\mH$ and every set $S \subseteq V(\mH)$, let $d_{\mH}(S)$ be the number of edges in $\mH$ containing $S$, and for every integer $i \in \{1,2,\ldots, r\}$, let $\Delta_i(\mH)$ be the maximum value of $d_\mH(S)$ over all sets $S \subseteq V(\mH)$ with $|S| = i$.
\begin{conjecture}[Balanced supersaturation] \label{conj::generalBalSup}
    For every integer $d \ge 2$, there is a positive real number $T = T(d)$ such that the following is true as $q$ goes to infinity.
    For every set $U \subseteq \F_q^d$ with $|U| \ge Tq$, there is a collection $\mH_U \subseteq \S_{d+1}(U)$ such that
    \begin{enumerate}[(i)]
        \item $|\mH_U| = \Omega(|U|^{d+1} / q)$,
        \item $\Delta_i(\mH_U) = O\left(|U|^{d+1-i} / q^{1- \frac{i-1}{d}}\right)$ for every integer $i\in \{  1, 2,\ldots,   d\}$,
    \end{enumerate}
    where all the implicit constants depend only on $d$.
\end{conjecture}

\begin{claim} \label{balSupImpMainConj}
    For every integer $d \ge 2$, if \cref{conj::generalBalSup} holds, then \cref{conj::generalTheorem} also holds.
\end{claim}

The main body of this paper is to prove  \cref{conj::generalBalSup}  for $d=3$, which, by \cref{balSupImpMainConj}, determines the order of magnitude of $\alFqdp[3]$ up to a polylogarithmic factor, see \cref{fig::alpha3}.  This also solves  Problem 6.1 in~\cite{chen2023random}. We remark that the connection between \cref{conj::generalBalSup} and \cref{conj::generalTheorem} is implicitly mentioned in~\cite{chen2023random} and one of their many results is that \cref{conj::generalBalSup} holds for $d=2$.
\begin{theorem} \label{3dBalSup}
    \cref{conj::generalBalSup} holds for $d=3$.
\end{theorem}

\begin{figure}[ht]
\centering

\begin{tikzpicture}
    \newcommand {\xl} {4}
    \newcommand {\xo} {0.03}
    \newcommand {\xt} {0.5}
    \newcommand {\xth} {2.5}
    \newcommand {\xf} {3.4}
    \newcommand {\yl} {11}
    \newcommand {\ys} {3.33}
    \newcommand {\yt} {3.22}
    \newcommand {\yth} {3.41}
    \newcommand {\yf} {10}

    \begin{axis}[
    width = 0.7\textwidth,
    height = 0.4\textwidth,
    xlabel = {$p$},
    ylabel = {$\alpha(\mathbb{F}^3_q, p)$},
    xmin = 0, xmax = \xl,
    ymin = 0, ymax = \yl,
    xtick = {\xo, \xt, \xth, \xf},
    xticklabels = {$q^{-3}$, $q^{-8/3}$, $q^{-2/3}\log^2 q$, $1$},
    ytick = {\ys, \yf},
    yticklabels = {$q^{1/3+o(1)}$, $q$},
    axis lines=middle,
    xlabel style={at=(current axis.right of origin), anchor=west},
    ylabel style={at=(current axis.above origin), anchor=south},
    ]
    \addplot[color = black]
        plot coordinates {
            (\xo,0)
            (\xt,\yt)
            (\xth,\yth)
            (\xf,\yf)
        };

    \addplot[color = gray,dashed]
        plot coordinates {
            (\xt, 0)
            (\xt, \yt)
        };
    \addplot[color = gray,dashed]
        plot coordinates {
            (0, \ys)
            (\xt, \ys)
        };
    \addplot[color = gray,dashed]
        plot coordinates {
            (\xth, 0)
            (\xth, \yth)
        };
    \addplot[color = gray,dashed]
        plot coordinates {
            (\xf, 0)
            (\xf, \yf)
        };
    \addplot[color = gray,dashed]
        plot coordinates {
            (0, \yf)
            (\xf, \yf)
        };
    \end{axis}

\end{tikzpicture}

\caption{The behavior of $\alFthrdp$ in terms of $p$.} 
\label{fig::alpha3}
\end{figure}

While, by now, there are balanced supersaturation theorems proved in various set-ups (see for example~\cite{morris2016number, ferber2020supersaturated, corsten2021balanced, mckinley2023random, jiang2024tree, nie2024random, Jiang2024Balanced}), somehow there is no unified ``one size fits all'' type of strategy. We trust that our probabilistic proof method of \cref{3dBalSup} could be adapted to attack other problems.

The rest of this paper is organized as follows. In \cref{sec::bal3d}, we prove \cref{3dBalSup}. In \cref{sec::balSupImpMainConj}, we prove \cref{balSupImpMainConj}. In \cref{sec::conRemark}, we give some concluding remarks about some barriers in proving \cref{conj::generalBalSup} for $d > 3$.

\textbf{Notation.}
For a hypergraph $\mH$, we use $v(\mH)$ for its number of vertices and $e(\mH)$ or $|\mH|$ for its number of hyperedges.
For a set $S$ and a non-negative integer $k$, we use $\b{S}{k}$ for the collection of all subsets of $S$ with size $k$.
We use the standard $O$, $\Omega$, $\Theta$, $o$, $\omega$ notation: $f(x) = O(g(x))$ and $g(x) = \Omega(f(x))$ if there is a constant $C > 0$ such that $f(x) \le Cg(x)$, $f(x) = \Theta(g(x))$ if both $f(x) = O(g(x))$ and $f(x) = \Omega(g(x))$, and $f(x) = o(g(x))$ and $g(x) =\omega(f(x))$ if $\lim_{x \to \infty} f(x)/g(x) = 0$.

\section{Balanced supersaturation lemma for \texorpdfstring{$d = 3$}{d=3}} \label{sec::bal3d}
In this section, we  prove \cref{3dBalSup}.
We remark that the same proof also works for $d=2$, and we will discuss some of the obstacles in higher dimensions in \cref{sec::conRemark}.
Our main idea for proving \cref{3dBalSup} is to analyze the structure of sets based on  at which ``dimension'' the degeneracy mainly occurs, namely, whether there are many triples in a line in $U$ or not, see \cref{manynkcaoverq}. In each case, we show the existence of the family $\mH_U$ satisfying the properties in \cref{conj::generalBalSup} by constructing it randomly.

We will use the following version of the Chernoff bound, see Corollary~A.1.14 in \cite{alon2016probabilistic}.

\begin{lemma} \label{chernoff}
    Let $Y$ be the sum of mutually independent indicator random variables, $\mu = \E[Y]$. For all $\eps > 0$,
    \[
        \P [|Y - \mu| > \eps \mu] < 2e^{-c_\eps \mu},
    \]
    where $c_\eps = \min \{ -\ln(e^{\eps} (1+\eps)^{-(1+\eps)}),\,\eps^2/2  \}$.
\end{lemma}
\noindent
Note that when $\eps > 10$ in \cref{chernoff}, we have $c_\eps > \eps$, hence $\P[ Y > (\eps+1) \mu] < 2e^{-\eps \mu}$. Therefore, for  $h \ge 20 \mu > 0$, we have
\begin{equation} \label{ourChernoff}
    \P[Y > h] = \P[Y > (h/\mu-1+1) \mu] < 2e^{-(h/\mu -1) \mu} = 2e^{-(h-\mu)} \le 2e^{-h/2}
    .
\end{equation}

\begin{proof}[Proof of \cref{3dBalSup}]
    Fix $d = 3$ and let $T$ be a sufficiently large real number. Let $U$ be an arbitrary subset of $\F_q^3$ with $|U| = n \ge Tq$. In particular, when we write in the proof $\Omega(n/q)$, we will always assume that it is much larger than $1$. For every set $S\subseteq \F_q^3$, we write $n_S \ce |U\cap S|$. Our aim is to prove that there is a family $\hu \subseteq \S_{4}(U)$ such that $|\hu| = \Omega(n^4/q)$,
    $\Delta_1(\hu) = O(n^3 / q)$, $\Delta_2(\hu) = O(n^2 / q^{2/3})$, and
    $\Delta_3(\hu) = O(n / q^{1/3})$.

    We use the following notation.
    For every integer $i \in \{1,2,3\}$ and every set $A \in \G_i(U)$, let $K_A$ be the $(i-1)$-flat containing $A$, $\mP_A$ be the set of  $i$-flats containing $K_A$, and
    \[
        \kca \ce \ka \sm \bigcup_{B \subsetneq A} \ka[B].
    \]

    Assume first that  $|\S_{4}(U)| = \Omega(n^4)$. Then, let $\hu$ be the set where every set in $\S_{4}(U)$ is included with probability $p = 1/q$, independently from each other. We have
    $\ehu = p|\S_{4}(U)| = \Omega(n^4/q)$ and $\E[d_{\hu}(S)] \le n^{4-|S|} \cdot 1/q = O(n^{4-|S|} / q^{1-(|S|-1)/3})$ for every set $S \subseteq U$ with size $|S| \in \{1,2,3\}$.
    By \cref{chernoff}, we have with high probability that $|\hu| = \Omega(n^4/q)$.
    Also, note that $n^{4-|S|} / q^{1-(|S|-1)/3} \ge q^{2/3}$ for every $S\subseteq U$ with $|S| \in \{1,2,3\}$ and there are at most $q^{O(1)}$ choices for all such sets $S$.
    By the union bound and \cref{ourChernoff}, we have that (ii) is satisfied with probability at least $1-q^{O(1)}e^{-\Omega(q^{2/3})}$. Therefore, there is a choice of $\hu \subseteq \S_4(U)$ satisfying both (i) and (ii).

    Hereinafter, we assume that $|\S_{4}(U)| = o(n^{4})$. In particular, for every $i$-flat $K$ ($1\le i \le 2$), we have
    \begin{equation} \label{equ::nk=on}
        n_K = o(n).
    \end{equation}

    \begin{claim} \label{manynkanvoerq}
        There is a family $\B \subseteq \G_{3}(U)$ with $|\B| = \Omega(n^3)$ such that $\nka[B] = \Omega(n/q)$ for every $B \in \B$.
    \end{claim}
    \begin{proof}
        For every pair $S  \in \b{U}{2}$, by \cref{equ::nk=on}, we have $n_{U\sm K_S} = n - n_{K_S} \ge n/2$. Note that
     $\{(U \cap P) \sm K_S: P \in \mP_S\}$ partitions
        $U \sm K_S$  into $q+1$ sets, so $\sum_{P \in \mP_S} n_{P \sm K_S} \ge n/2$ and
        \begin{equation} \label{ave}
            \sum_{\substack{P \in \mP_S \\ n_{P \sm K_S} = \Omega(n/q)}} n_{P \sm K_S} = \Omega(n).
        \end{equation}
        We let
        \[
            \B \ce \left\{S \cup \{x\}: S\in \b{U}{2},\,x\in (U \cap P) \sm K_S \textrm{ where }P\in \mP_S \textrm{ and } n_{P \sm K_S} = \Omega(n/q)\right\}.
        \]
        Note that each set in $\B$ is generated at most three times. By~\eqref{ave},
        we have $|\B| = \Omega(n^3)$.
    \end{proof}

    \begin{claim} \label{manynkcaoverq}
        There is an integer $j \in \{2,3\}$ such that the following is true. There is a family $\A_j \subseteq \G_{j}(U)$ with $|\A_j| = \Omega(n^j)$ such that for every $A \in \A_j$, we have
        \begin{equation} \label{equ::nkcanca}
            \nkca \ge \nka / 10 = \Omega(n/q).
        \end{equation}
    \end{claim}
    \begin{proof}
        Take the family $\B$ from \cref{manynkanvoerq}. For every $B = \{x,y,z\} \in \B$,  the $2$-flat $K_B$ is partitioned into $7$ sets:
        \[
            \{x\}, \{y\}, \{z\}, K^\circ_{x,y}, K^\circ_{x,z}, K^\circ_{y,z}, K^\circ_{B}.
        \]
        Let $f(B)$ be the set with the maximum number of points in $U$ among these $7$ sets. Since $\nka[B]= \Omega(n/q)$ and $T$ is sufficiently large, by choosing the implicit constant properly, we have that $f(B) \in \{K^\circ_{x,y}, K^\circ_{x,z}, K^\circ_{y,z}, K^\circ_{B}\}$.
        If there are $\Omega(n^3)$ sets $B \in \B$ such that $f(B) = \kca[B]$, then we are done with $j = 3$ and $\A_3$ being the collection of all such sets $B$. Otherwise, note that every pair is  in at most $n$ triples from $\B$, hence we are done with $j = 2$.
    \end{proof}

    According to the value of $j$ in \cref{manynkcaoverq}, we divide the remaining proof into the following cases. In each case, we give our probabilistic construction for $\hu$ and first prove that $\ehu = \Omega(n/q)$ and $\edhus = O(n^{4-|S|} / q^{1-(|S|-1)/3})$ for every $S \subseteq U$ with size $S \in \{1,2,3\}$. At the end of this proof, we show that a simple application of the union bound and the Chernoff bound guarantees the existence of $\hu$ satisfying (i) and (ii) in all cases.

    \textbf{Case 1:} $j = 3$ in \cref{manynkcaoverq}.

    In this case, we let $\A \ce \A_3$ from \cref{manynkcaoverq}. We build our $4$-sets in $\hu$ by adding a point in $\kca$ to $A$ for every $A \in \A$.

    For every $A \in \A$, let $\mF_{A} \ce \{A \cup \{x\}: x \in \kca \cap U\}$ and $\ha$ be the set where every set in $\mF_A$ is included with probability $p_A = \Theta(n/(q \nkca))\le 1$, independently of each other. Let $\hu \ce \cup_{A \in \A} \ha$.
     Note that $\hu \subseteq \S_4(U)$, since $x \in \kca \cap U$ for every $A$.
     For every $A \in \A$, we have
     \[
        \eha = |\mF_A| \cdot p_A = \nkca \cdot p_A = \Omega\left(\frac{n}{q} \right),
     \]
     hence
     \[
        \ehu \ge \frac{1}{4} \sum_{A \in \A} \eha = |\A| \cdot  \Omenoverq = \Omenpoverq{4}.
     \]

     For (ii),
     it suffices to prove that $\edhus = O(n/q)$ for every $S \subseteq U$ with $|S| = 3$.
     Indeed, it implies $\edhus[S'] \le n \cdot O(n/q) = O(n^2/q) = O(n^2 /q^{2/3})$ for every $S' \subseteq U$ with $|S'| = 2$, since every pair $S'$ is trivially in at most $n$ triples in $\b{U}{3}$, and similarly $\edhus[S'] \le n^2 \cdot O(n/q) = O(n^3/q)$ for every $S' \subseteq U$ with $|S'| = 1$.
     Fix a set $S \subseteq U$ with $|S| = 3$. If $S \in \S_3(U)$, then $\dhus = 0$ by the definition of $\mF_A$. If $S \in \G_3(U)$, then let $\A' \subseteq \A$ be the collection of sets $A \in \A$ such that $K_A = K_S$ and $|A \sm S| = 1$.
     Note that $|\A'| \le 3\nka[S]$ and if $A \in \A \sm (\A' \cup \{S\})$, then $\dfas = 0$. By the definition of $\A$, we have
     \[
         \sum_{A \in \A'} \edhas
         = \sum_{A \in \A'} p_{A}
         = \sum_{A \in \A'} O\left(\frac{n}{q\nkca}\right)
         = \sum_{A \in \A'} O\left(\frac{n}{q\nka[S]}\right)
        \le 3n_{K_S} \cdot O\left(\frac{n}{q\nka[S]}\right)
        = \Onoverq
        .
     \]
     If $S \notin \A$, then we are done. Otherwise, we still have
     \[
        \edhus \le \E[d_{\mH_S}(S)] +
        \sum_{A \in \A'} \edhas
        = \nkca[S] \cdot p_S + O\left(\frac{n}{q}\right) = O\left(\frac{n}{q}\right).
     \]

    \textbf{Case 2:} $j = 2$ in \cref{manynkcaoverq}.

    Let $\A_2$ be the set from \cref{manynkcaoverq}.
    Note that we can simply assume that $\A_2 = \{A \in \b{U}{2}: n_{K_A} = \Omega(n/q)\}$.
    We further divide the proof into two cases, depending on whether there is a subset $\A'_2 \subseteq \A_2$ with size $\Omega(n^2)$ such that every $A \in \A'_2$ has $\nka \ge n/q^{2/3}$.

    \textbf{Subcase 2.1:} Such $\A_2'$ does not exist.

    Then, let $\A \ce \{A \in \A_2: \nka < U/q^{2/3}\}$. Note that we have $|\A| = (1-o(1))|\A_2| = \Omega(n^2)$ and $\nka = \Omega(n/q) \ge 3$ for every $A \in \A$ by our assumption on $\A_2$.
    In this case,
    for every $A \in \A$, we add two additional points from every plane containing $A$ while we require that no three points are in a line in order to guarantee (ii).
    Note that for every $A \in \A$ and every pair of points $\{x,y\}$ which is disjoint from $K_A$ but coplanar with $A$, there is at most one point in $K_A$ which is in the same line with $x,y$. By the fact that $\nka \ge 3$, we are still able to find at least one-third of the pairs of points from $K_A$ to form a $4$-tuple eligible to include in $\hu$.

    For every $A \in \A$ and $P \in \mP_A$ with $n_{P \sm K_A} = \Omega( n/q)$, let
    \[
    \mF_{A,P} \ce \{F : |F| = 4,\,A \subseteq F \subseteq P,\, \textrm{ every subset of $F$ of  size $3$ is in $\G_3(U)$}\}
    \]
    and $\mH_{A,P}$ be the set where every set in $\mF_{A,P}$ is included with probability $p_{A,P} \ce \Theta(n/(qn_{P\sm K_A}))$, independently of each other. Let
    \[
        \mH_A \ce \bigcup_{\substack{P: P\in \mP_A\\
        n_{P\sm A} = \Omega(n / q)
        }} \mH_{A,P}
        \quad\quad \quad\textrm{and}\quad \quad\quad
        \mH_U \ce
        \bigcup_{A:A \in \A}
        \mH_A
        .
    \]

    To prove (i),
    recall $3 \le \nka < n/q^{2/3}$ for every $A \in \A$. Then we have
    \[
        \sum_{\substack{K: \textrm{$K$ is a $1$-flat} \\ 3 \le n_K < n/q^{2/3}}} \b{n_K}{2}
        \ge |\A| = \Omega(n^2).
    \]
    For every $1$-flat $K$ with $3 \le n_K < n/q^{2/3}$, a point $z \in K \cap U$, and every $2$-flat $P$ with $K \subseteq P$, we let
    \[
        M_{z, K, P} \ce \left\{ \{a,b\} \in \b{(P \sm K) \cap U}{2}:
        z \in K_{\{a,b\}}
        \right\} ,
    \]
    the pairs $\{a,b\}$ of points in $(P \sm K) \cap U$ such that $\{z,a,b\} \in \S_3(U)$,
    and $m_{z, K, P} \ce |M_{z, K, P}|$.
    Let $z_{K,P}$ be the point in $K \cap U$ with the largest $m_{z,K,P}$. Note that for $z_1 \neq z_2 \in K$, we have $M_{z_1,K,P} \cap M_{z_2,K,P} = \es$.
    Therefore, for every pair $\{x,y\}\subseteq K$ where $z_{K,P} \notin \{x,y\}$, we have
    \begin{equation} \label{mxKP}
        m_{x,K,P} + m_{y,K,P} + m_{z_{K,P}, K,P} \le \b{n_{P\sm K}}{2}
        \quad \textrm{ implying} \quad
        m_{x,K,P} + m_{y,K,P} \le \frac{2}{3} \b{n_{P\sm K}}{2}.
    \end{equation}
    Also, by the assumption that $n_K \ge 3$, we have that for every $2$-flat $P \supseteq K$, there are
    \begin{equation} \label{nK-12nK2}
        \b{n_K-1}{2} = \frac{\b{n_K-1}{2}}{\b{n_K}{2}}\cdot \b{n_K}{2} = \frac{n_K-2}{n_K} \cdot \b{n_K}{2} \ge \frac{1}{3} \b{n_K}{2}
    \end{equation}
    pairs $\{x,y\} \subseteq K \cap U$ such that $z_{K,P} \notin \{x,y\}$.
    Now, for every $1$-flat $K$ with $3 \le n_K < n/q^{2/3}$, we have
    \begin{align*}
        \sum_{A \subseteq K} \eha
        &= \Omega\Bigg(
        \sum_{\substack{P: \textrm{$P$ is a $2$-flat}\\ K \subseteq P \\n_{P\sm K} = \Omega(n/q) }}\ \
        \sum_{\substack{\{x,y\}: \{x,y\} \subseteq K \cap U\\ z_{K,P} \notin \{x,y\} }}
        \Bigg( \b{n_{P \sm K}}{2} - m_{x,K,P} - m_{y,K,P} \Bigg) \cdot p_{K,P}
        \Bigg)\\
        & \stackrel{\eqref{mxKP}}{=} \Omega \Bigg(
        \sum_{\substack{P: \textrm{$P$ is a $2$-flat}\\ K \subseteq P\\n_{P\sm K} = \Omega(n/q)} }
        \ \ \sum_{\substack{\{x,y\}: \{x,y\} \subseteq K \cap U\\ z_{K,P} \notin \{x,y\} }}
        \b{n_{P \sm K}}{2}  \cdot \frac{n}{q n_{P \sm K}}
        \Bigg)\\
        &\stackrel{\eqref{nK-12nK2}}{=} \Omega \Bigg(
        \sum_{\substack{P: \textrm{$K$ is a $2$-flat}\\ K \subseteq P\\n_{P\sm K} = \Omega(n/q)} }
        \b{n_K}{2}
        \cdot \b{n_{P \sm K}}{2}  \cdot \frac{n}{q n_{P \sm K}}
        \Bigg) \stackrel{\eqref{equ::nk=on}}{=} \Omega \Bigg( \b{n_K}{2} \cdot n \cdot \frac{n}{q} \Bigg)
        .
    \end{align*}
    Each $4$-tuple could be included at most $\binom{4}{2}$ times, therefore,
    \[
        \ehu \ge \frac{1}{\b{4}{2}} \sum_{A \in \A} \eha = \Omega\left(n^2 \cdot n \cdot \frac{n}{q}\right) = \Omenpoverq{4}.
    \]

    For (ii), fix an arbitrary set $S \subseteq U$ with $|S| = i \in \{1,2,3\}$. For $i= 3$, note that if $S \in \S_3(U)$, then $\dhus = 0$ by the definition of $\fa$, so we can assume that $S = \{x,y,z\} \in \G_3(U)$. Let $P = K_S$.
    Then, $\dfas > 0$ can only happen when $A \subseteq P$ and $|A \sm S| \le 1$. If $|A \sm S| = 0$, then $A = \{x,y\}$ or $\{x,z\}$ or $\{y,z\}$. For $A = \{x,y\}$, we have $d_{\mF_{A,P}}(S) = n_{K^\circ_S} \le n_{P \sm K_A}$ and hence
    \[
        \edhas =
        \E[d_{\mH_{A,P}}(S)] = d_{\mF_{A,P}}(S) \cdot p_{A,P} = O\left( n_{P \sm K_A} \cdot \frac{n}{qn_{P \sm K_A}} \right)
        = \Onoverq.
    \]
    The same arguments hold for $A = \{x,z\}$ or $A = \{y,z\}$. If $|A \sm S|  = 1$, then we first consider those sets $A$ containing $x$ but not $y$ or $z$.
    Let $\K$ be the collection of $1$-flats $K \subseteq P$ with $x \in K \subseteq P$, $y,z \notin K$, $n_K \le n/q^{2/3}$, and $n_{P \sm K} = \Omega(n/q)$.
    By our construction, we have
    \begin{equation} \label{equ::ranDelta3}
        \sum_{\substack{A \in \A \\ x\in A, \ y,z \notin A}} \edhas \le
        \sum_{K \in \K}
        n_{K \sm \{x\}} \cdot \Theta\left( \frac{n}{qn_{P \sm K}} \right)
        =
        O \Bigg(\frac{n}{q} \cdot \sum_{K \in \K}  \frac{n_{K \sm \{x\}}}{ n_{P\sm K}} \Bigg)
        .
    \end{equation}
    Assume that $\hat{K}$ is the $1$-flat in $\K$ which has the largest $n_K$. For $\hat{K}$, it contributes to~\eqref{equ::ranDelta3}  at most $O(n_{\hat{K}\sm\{x\}}) \le O(n/q^{2/3})$, since $n_{P \sm K} = \Omega(n/q)$.
     For all other $K \neq \hat{K} \in \K$,  since $n_{P \sm K} \ge n_{\hat{K} \sm \{x\}} \ge n_{K \sm \{x\}}$, we have
     \[
        \frac{n_{K \sm \{x\}}}{n_{P \sm K}}
        \le \frac{2n_{K \sm \{x\}}}{n_{P \sm K}+n_{K \sm \{x\}}} = \frac{2n_{K \sm \{x\}}}{n_{P \sm \{x\}}},
     \]
     hence
     \[
     \sum_{\substack{K\in \K \\ K \neq \hat{K}}}
        \frac{n_{K \sm \{x\}}}{n_{P \sm K}}
        \le
        \sum_{\substack{K\in \K \\ K \neq \hat{K}}} \frac{2n_{K \sm \{x\}}}{n_{P \sm \{x\}}}
        \le 2 \cdot \frac{1}{n_{P \sm \{x\}}} \cdot \sum_{K \in \K} n_{K \sm \{x\}}
        \le 2 \cdot \frac{1}{n_{P \sm \{x\}}} \cdot n_{P \sm \{x\}}
        \le 2
        .
     \]
     Continuing \cref{equ::ranDelta3}, we have
     \[
        \sum_{\substack{A \in \A \\ x\in A, \ y,z \notin A}} \edhas = O\left( \frac{n}{q^{2/3}} \right) + O\left( \frac{n}{q} \right) = O\left( \frac{n}{q^{2/3}} \right).
    \]
    The same arguments hold for the sets $A \in \A$ containing only $y$ or $z$, so $\sum_{A \in \A, \ |A \sm S| = 1}\edhas =  O(n/q^{2/3})$. Therefore
    \begin{equation} \label{equ::case21i3}
    \edhus \le \sum_{A \in \A} \edhas = O\left( \frac{n}{q} \right) + O \left( \frac{n}{q^{2/3}} \right) = O \left( \frac{n}{q^{2/3}} \right),
    \end{equation}
    which is $O(n/q^{1/3})$, as claimed.

    For $i = 2$, we  use \cref{equ::case21i3} and get $\edhus \le n \cdot O(n/q^{2/3}) = O(n^2/q^{2/3})$.

    For $i = 1$, we have two cases, either $A \supseteq S$ or not. There are  at most $n$ sets $A \in \A$ containing $S$, and for each of them, we have
    \begin{align*}
        \edhas = \eha
        &= \sum_{\substack{P\in \mP_A \\n_{P\sm K_A} = \Omega(n/q) }} \E[|\mH_{A,P}|] \le
        \sum_{\substack{P\in \mP_A \\n_{P\sm K_A} = \Omega(n/q) }}
        \b{n_{P \sm K_A}}{2}  \cdot p_{A,P} \\
        &= O \Bigg( \frac{n}{q} \cdot \sum_{\substack{P\in \mP_A \\n_{P\sm K_A} = \Omega(n/q) }} n_{P \sm K_A} \Bigg)
        = O \left( \frac{n^2}{q} \right)
        .
    \end{align*}
    There are  at most $n^2$ sets $A \in \A$ not containing $S$. If $S \subseteq K_A$, then $\dfas = 0$. If $S \not\subset K_A$, then let $P = K_{A \cup S}$. Note that $\dfas > 0$ only if $n_{P \sm K_A} = \Omega(n/q)$, due to our definition of $\ha$, and then we have
    \[
        \edhas = \E[d_{\mH_{A,P}}(S)] = n_{K^\circ_{A \cup S}} \cdot \Theta\left(\frac{n}{qn_{P\sm K_A}}\right) \le n_{P\sm K_A} \cdot \Theta\left(\frac{n}{qn_{P\sm K_A}}\right) = O\left( \frac{n}{q} \right).
    \]
    Therefore, we have
    \[
        \edhus \le \sum_{A \in \A} \edhas \le n\cdot O\left( \frac{n^2}{q} \right) + n^2 \cdot O\left(\frac{n}{q}\right) = O\left( \frac{n^3}{q} \right),
    \]
    as claimed.

    \textbf{Subcase 2.2:} Such $\A_2'$ exists.

    Let $\A \ce \A_2' = \{A \in \A_2: \nka > |U|/q^{2/3}\}$, so $|\A| = \Omega(n^2)$ by assumption.
    In this case, for every $A \in \A$, we add to it one additional point in $K_A$ and one additional point in $U \sm K_A$.

    For every $A \in \A$, we let
    \[
        \mF_{A} \ce \{A \cup \{x\} \cup \{y\}: x \in \kca,\, y\in U \sm \ka\}
    \]
    and $\ha$ be the set where each set in $\fa$ is included with probability $p_A = \Theta(n/(q\nkca))$,  independently from each other. Let $\hu \ce \cup_{A \in \A} \ha$.

    To prove (i), observe that for every $A \in \A$, by \cref{equ::nk=on}, we have
    \[
        \eha = |\fa|\cdot p_A = \nkca \cdot |U \sm K_A| \cdot \Theta\left(\frac{n}{q\nkca}\right) =  \Omega(n) \cdot \Omenoverq  = \Omenpoverq{2},
    \]
    and
    \[
        \ehu \ge \frac{1}{\b{3}{2}} \cdot \sum_{A \in \A} \eha
        = \Omega(n^2) \cdot \Omenpoverq{2} = \Omenpoverq{4}.
    \]

    To prove (ii), fix an arbitrary set $S \subseteq U$ with $|S| = i$ ($1\le i\le 3$). First, we consider the case where $i=3$. If $S \in \S_3(U)$, then $\dfas > 0$ only if $A \subseteq S$; there are at most three such sets $A$. For each of them, we trivially have $\dfas \le n$, implying  $\edhas \le n \cdot p_{A} = O(n/q^{1/3})$ due to our assumption that $\nka \ge n/q^{2/3}$ and our definition of $p_A$.
    Hence, we have $\edhus = O(n/q^{1/3})$. If $S =\{x,y,z\} \in \G_3(U)$, then $\dfas > 0$ only if $A \subseteq K_{\{x,y\}}$, $K_{\{x,z\}}$ or $K_{\{y,z\}}$. For those sets $A \subseteq K_{\{x,y\}}$, $A$ can be $\{x,y\}$, $\{x,w\}$, or $\{y,w\}$, where $w \in \kca[\{x,y\}]$. For all the sets $A$ of the form $\{x,w\}$ or $\{y,w\}$, their contribution to $\edhus$ is at most
    \[
        \sum_{w \in K_{\{x,y\}}^\circ} \Theta\left( \frac{n}{q\nkca[\{x,w\}]} \right) + \Theta\left( \frac{n}{q\nkca[\{y,w\}]} \right)
        \le 2\nka[\{x,y\}] \cdot O\left( \frac{n}{q\nka[\{x,y\}]}\right)
        = O\left( \frac{n}{q} \right).
    \]
    The contribution of the set $A = \{x,y\}$  to $\edhus$ is at most
    \[
        \nkca[\{x,y\}] \cdot p_{\{x,y\}}
        = \nkca[\{x,y\}] \cdot \Theta\left( \frac{n}{q\nkca[\{x,y\}]}\right) = O\left( \frac{n}{q}\right).
    \]
    The same arguments hold for the sets $A$ from $K_{\{x,z\}}$ and $K_{\{y,z\}}$.
    Hence, we have
    \begin{equation} \label{equ::Sg3Unoverq}
        \edhus = O(n/q) \textrm{ for every } S \in \G_3(U).
    \end{equation}
    Therefore, we have $\edhus = O(n/q^{1/3})$ for every $S \subseteq U$ with $|S| = 3$.

    For $i = 2$, note that all the hyperedges in $\hu$ containing $S$ must contain some vertex $z$ which is not in $K_S$. By \cref{equ::Sg3Unoverq}, we have that $\edhus[S] = n \cdot O(n /q) = O(n^2/q)$.

    For $i=1$, we use the result for $i=2$ and have $\edhus \le n \cdot O(n^2 /q) = O(n^3 /q)$.

    The remaining part of the proof is to prove the actual existence of $\hu \subseteq \S_4(U)$ satisfying both (i) and (ii) in all cases. In the discussion above, we proved that in each case, there is a probabilistic construction for $\hu$ such that $\ehu = \Omega(n^4/q)$ and $\edhus = O(n^{4-|S|} / q^{1-(|S|-1)/3})$ for every $S \subseteq U$. Then, a simple application of the Chernoff bound guarantees the actual existence of such $\hu$. Indeed, in both Cases 1 and 2, we have $\hu = \cup_{A \in \A} \ha$. Note that for every $A \in \A$ and $S\subseteq U$, we have that $|\ha|$ and $\dhas$ are the sum of mutually independent indicator random variables. This is easy to see in Case~1 and Case 2.2, as $\ha$ is a subset of $\fa$ where each element is kept with a certain probability independently of each other. For Case 2.1, it suffices to notice that $\mF_{A, P_1} \cap \mF_{A, P_2} = \es$ for $P_1 \neq P_2 \in \mP_A$ and then recall our definition of $\ha$. Therefore, $\sum_{A \in \A} |\mH_A|$ and $\sum_{A \in \A} \dhas$ are also sum of mutually independent indicator random variables.
    We also have
    \begin{equation} \label{equ::husumAhaO1hu}
        |\hu| \le \sum_{A \in \A} |\ha| \le O(1) \cdot |\hu|
    \end{equation}
    and
    \begin{equation} \label{equ::dhussumAdhasO1dhus}
        \dhus \le \sum_{A \in \A} \dhas \le O(1) \cdot \dhus,
    \end{equation}
    since every hyperedge in $\hu$ is trivially counted at most $O(1)$ times due to our construction. Hence, we have
    $\sum_{A \in \A}\eha = \Omega(n^4/q)$ and $\sum_{A \in \A} \edhas = O(n^{4-|S|} / q^{1-(|S|-1)/3})$ for every $S \subseteq U$.

    The remaining calculation is the same as the case at the beginning of the proof where $|\S_{4}(U)| = \Omega(n^4/q)$.
    By \cref{chernoff}, we have $\sum_{A \in \A} |\mH_A| = \Omega(n^4/q)$ with high probability, and by \cref{equ::husumAhaO1hu}, so is $|\hu|$. We have $n^{4-|S|} / q^{1-(|S|-1)/3} \ge q^{2/3}$ for every $S\subseteq U$ with $|S| \in \{1,2,3\}$ and there are at most $q^{O(1)}$ choices for all such sets $S$.
    By the union bound and \cref{ourChernoff}, with probability at least $1-q^{O(1)}e^{-\Omega(q^{2/3})}$ we have $\dhus \le \sum_{A \in \A} \dhas = O(n^{4-|S|} / q^{1-(|S|-1)/3})$ for every $S \subseteq U$ with $|S| \in \{1,2,3\}$.
    Thus, there is a choice of $\hu \subseteq \S_4(U)$ satisfying both (i) and (ii). \qedhere

    \end{proof}

\section{Proof of \texorpdfstring{\cref{balSupImpMainConj}}{Claim 3}} \label{sec::balSupImpMainConj}
The proof method used in this section is by now a standard application of the Hypergraph Container Method~\cite{balogh2015independent, saxton2015hypergraph}.
We use the following version of the Hypergraph Container Lemma, see Corollary~2.8 in~\cite{janzer2025improved} for its proof.
\begin{lemma}
\label{lem::hypContainerLemma}
For every positive integer $s$ and positive real numbers $x$ and $\lambda$, the following holds. Suppose that $\G$ is an $s$-uniform hypergraph such that $x v(\G)$ and $v(\G) / \lambda$ are integers, and for every $\l \in \{1,2,\ldots, s\}$,
\begin{equation} \label{equ::codContainerLemma}
\Delta_{\ell}(\G) \leq \lambda \cdot x^{\ell-1} \frac{e(\G)}{v(\G)}.
\end{equation}
Then, there exists a collection $\C$ of at most $v(\G)^{s\cdot x\cdot v(\G)}$ sets of size at most $\left(1-\delta \lambda^{-1}\right) v(\G)$ such that for every independent set $I$ in $\G$, there exists some $R \in \C$ with $I \subseteq R$, where $\delta=2^{-s(s+1)}$.
\end{lemma}

The main idea for \cref{conj::generalTheorem} is to apply \cref{lem::hypContainerLemma} to the collections $\mH_U$ in \cref{conj::generalBalSup}. The property (ii) in \cref{conj::generalBalSup} guarantees that \cref{equ::codContainerLemma} is always satisfied. We will apply \cref{lem::hypContainerLemma} for a finite number of times and finally show that there is a small collection $\C$ of sets with size at most $T(d) \cdot d$ such that every set in general position in $\F_q^d$ is contained in some set in $\C$. This enables us to give an upper bound on the number of sets in general position in $\F_q^d$ and then control the maximum size of it in the random subspace.

\begin{proof}[Proof of \cref{conj::generalTheorem} assuming \cref{conj::generalBalSup}]

Fix a positive integer $d \ge 2$ and assume that \cref{conj::generalBalSup} holds for it. Let $T$ be the constant in \cref{conj::generalBalSup} and
let $C$ be sufficiently large. We first have the following claims.

\begin{claim} \label{largepO}
If $p \ge Cq^{-1+1/d} \log^2 q$, then $\alFqdp = O(pq)$ with high probability.
\end{claim}
\begin{proof}
    We first prove that for every integer $j \ge 0$, there is a collection $\C_j$ of at most  $(q^d)^{2(d+1)jq^{1/d}}$
    sets of size at most $\max(Tq, (1- \Omega(1))^j q^d)$
    such that every set $I\subseteq \F_q^d$ in general position is in some set from $\C_j$.
    Note that the base case $j = 0$ is trivial by taking $\C_0 = \{\F_q^d\}$.
    Given $\C_j$, we define $\C_{j+1}$ as follows. For every $U \in \C_j$ with $|U| > Tq$, 
    let $\hu$ be the collection given in \cref{conj::generalBalSup}. Note that if $I$ is a subset of $U$ in general position, then it is an independent set in $\hu$. We apply \cref{lem::hypContainerLemma} by letting $\G = \hu$, $\lambda = O(1)$, and $x$ be a real number in $[q^{1/d}/|U|, 2q^{1/d}/|U|]$ such that $xv(\G)$ and $v(\G)/\lambda$ are integers, and obtain  a collection $\C_U$. Let
    \[
        \C_{j+1} \ce \left\{U \in \C_j: |U| \le Tq \right\} \cup \bigcup_{U \in \C_j: |U| > Tq} \C_U.
    \]
    We have
    \[
        |\C_{j+1}| \le |\C_j| \cdot |U|^{(d+1)\cdot \frac{2q^{1/d}}{|U|} \cdot |U|} \le (q^d)^{2(d+1)j\cdot q^{1/d}} \cdot (q^d)^{2(d+1)\cdot q^{1/d}}
        = (q^d)^{2(d+1)(j+1)\cdot q^{1/d}}.
    \]
    For every $R \in \C_{j+1}$, we have either $|R| < Tq$ or there is a set $U \in \C_j$ such that $R \in \C_U$. In the latter case, we have
    \[
        |R| \le (1-\Omega(1)) |U| \le (1-\Omega(1)) (1-\Omega(1))^j q^d = (1-\Omega(1))^{j+1} q^d.
    \]
    Hence, $|R| \le \max(Tq, (1- \Omega(1))^{j+1} q^d)$.

    Now, let $j$ be a suitable value of order $\Theta(\log q)$ such that $(1-\Omega(1))^j q^d <Tq$. Let $t = \lceil Cpq \rceil \ge Cq^{1/d} \log ^2 q$. By the discussion in the last paragraph, the expected number of sets $I$ in general position with size $t$ in a $p$-random subset of $\F_q^d$ is at most
    \[
        q^{O(q^{1/d}\log q)} \b{Tq}{t} \cdot p^t
        \le q^{O(q^{1/d}\log q)} \left(\frac{eTqp}{t} \right)^t
        \le \left( \frac{1}{2} \right)^t =o(1).
    \]
    Therefore, by Markov's Inequality, we have that with high probability, there is no such $I$ in a $p$-random subset of $\F_q^d$, so $\alFqdp = O(pq)$.
\end{proof}

\begin{claim} \label{largepOmega}
If $p = \omega(q^{-1})$, then $\alFqdp = \Omega(pq)$ with high probability.
\end{claim}

\begin{proof}
    Recall the result of Erd\H{o}s that $\alFqd \ge q$ mentioned in \cref{sec::int}.
    We can fix a set $S\subseteq \F_q^d$ with $|S| = q$ which is in general position. \cref{chernoff} gives that with high probability, there are $\Omega(p|S|)$ elements of $S$ remaining in the $p$-random subset of $\F_q^d$, so $\alFqdp \ge \Omega(p|S|) = \Omega(pq)$.
\end{proof}

\begin{claim} \label{smallpO}
If $p = \omega(q^{-d})$, then $\alFqdp = O(pq^d)$ with high probability.
\end{claim}
\begin{proof}
Note that $|\F_q^d| = q^d$. By the assumption $p = \omega(q^{-d})$,
\cref{chernoff} gives that with high probability, there are $O(pq^d)$ points in the $p$-random subset of $\F_q^d$, so $\alFqdp = O(pq^d)$.
\end{proof}

\begin{claim} \label{smallpOmega}
If $\omega(q^{-d}) = p = o(q^{-d + 1/d})$, then $\alFqdp = \Omega(pq^d)$ with high probability.
\end{claim}
\begin{proof}
    By the assumption that $p = \omega(q^{-d})$ and \cref{chernoff}, we have with high probability that there are at least $pq^d/2$ points  in a $p$-random subset of $\F_q^d$. Note that $|\S_{d+1}(\F_q^d)| \le q^{d^2+d-1}$, since the  number of $(d-1)$-flats in $\F_q^d$ is $q^{d}$ and every $(d-1)$-flat contains at most $(q^{d-1})^{d+1} = q^{d^2-1}$ subsets of size $d+1$. Since $p = o(q^{-d + 1/d})$, we have $p^{d+1}q^{d^2+d-1} = o(pq^d)$, so by  Markov's Inequality, we have with high probability that there are at most $pq^d /4$ sets in $\S_{d+1}(\F_q^d)$ remaining in the $p$-random subset of $\F_q^d$, so
    \[
        \alFqdp \ge pq^d/2 - pq^d/4 = \Omega(pq^d). \qedhere
    \]
\end{proof}

Now, if $\omega(q^{-d}) = p = o(q^{-d+1/d})$, then by \cref{smallpO,smallpOmega}, we have $\alFqdp = \Theta(pq^d)$ with high probability. If $Cq^{-1+1/d}\log^2 q \le q \le 1$, then by \cref{largepO,largepOmega}, we have $\alFqdp = \Theta(pq)$ with high probability.
The last case where $q^{-d+1/d}/\log q \le p < Cq^{-1+1/d} \log^2 q$ simply follows from \cref{largepO}, \cref{smallpOmega}, and  monotonicity.
\end{proof}

\section{Concluding remarks} \label{sec::conRemark}
In this paper, we prove \cref{conj::generalBalSup} for $d = 3$. Unfortunately, we need new ideas for higher dimensions. Say, consider the case $d = 4$. We still have similar arguments as \cref{manynkanvoerq,manynkcaoverq}, which give that there is an integer $j \in\{ 2,3,4\}$ such that the conclusion of \cref{manynkcaoverq} holds. For the case $j = 3$ or $j = 4$, we can handle them by similar arguments as in the proof of \cref{3dBalSup}. For the case $j =2$, we now need to add three points from $U$ which are in the same $3$-flat with a given pair from $\A_2$.
However, it may happen that almost all such triples from the same $3$-flat are actually in a line, in which case we cannot guarantee (i) and (ii) at the same time. Intuitively, we think this case should not exist at all, but we do not know how to prove it. The same barriers exist for all higher dimensions.

\section*{Acknowledgment}
The first author is grateful for the hospitality of Hehui Wu at Fudan University (China), where he heard \cref{conj::generalBalSup} from Jiaxi Nie. The authors thank Jiaxi Nie for useful conversation on the project. We are grateful for the anonymous referees for their useful comments on the manuscript.


\begin{thebibliography}{10}

\bibitem{alon2016probabilistic}
N.~Alon and J.~H. Spencer.
\newblock {\em The probabilistic method}.
\newblock Wiley Series in Discrete Mathematics and Optimization. John Wiley \&
  Sons, Inc., Hoboken, NJ, fourth edition, 2016.

\bibitem{balogh2015independent}
J.~Balogh, R.~Morris, and W.~Samotij.
\newblock Independent sets in hypergraphs.
\newblock {\em J. Amer. Math. Soc.}, 28(3):669--709, 2015.

\bibitem{balogh2018number}
J.~Balogh and J.~Solymosi.
\newblock On the number of points in general position in the plane.
\newblock {\em Discrete Anal.}, Paper No. 16, 20, 2018.

\bibitem{brass2005research}
P.~Brass, W.~Moser, and J.~Pach.
\newblock {\em Research problems in discrete geometry}.
\newblock Springer, New York, 2005.

\bibitem{cardinal2017general}
J.~Cardinal, C.~D. T\'{o}th, and D.~R. Wood.
\newblock General position subsets and independent hyperplanes in {$d$}-space.
\newblock {\em J. Geom.}, 108(1):33--43, 2017.

\bibitem{chen2023random}
Y.~Chen, X.~Liu, J.~Nie, and J.~Zeng.
\newblock Random {T}ur{\'a}n and counting results for general position sets
  over finite fields.
\newblock {\em arXiv:2309.07744}, 2023.

\bibitem{corsten2021balanced}
J.~Corsten and T.~Tran.
\newblock Balanced supersaturation for some degenerate hypergraphs.
\newblock {\em J. Graph Theory}, 97(4):600--623, 2021.

\bibitem{dudeney1959amusements}
H.~E. Dudeney.
\newblock {\em Amusements in mathematics}.
\newblock Dover Publications, Inc., New York, 1959.

\bibitem{erdos1986some}
P.~Erd\H{o}s.
\newblock On some metric and combinatorial geometric problems.
\newblock {\em Discrete Math.}, 60:147--153, 1986.

\bibitem{ferber2020supersaturated}
A.~Ferber, G.~McKinley, and W.~Samotij.
\newblock Supersaturated sparse graphs and hypergraphs.
\newblock {\em Int. Math. Res. Not. IMRN}, (2):378--402, 2020.

\bibitem{furedi1991maximal}
Z.~F\"{u}redi.
\newblock Maximal independent subsets in {S}teiner systems and in planar sets.
\newblock {\em SIAM J. Discrete Math.}, 4(2):196--199, 1991.

\bibitem{janzer2025improved}
O.~Janzer and B.~Sudakov.
\newblock Improved {B}ounds for the {E}rd{\H{o}}s--{R}ogers
  ($s,s+2$)-{P}roblem.
\newblock {\em Random Structures \& Algorithms}, 66(1):Paper No. e21280, 2025.

\bibitem{Jiang2024Balanced}
T.~Jiang and S.~Longbrake.
\newblock Balanced supersaturation and {Tur}{\' a}n numbers in random graphs.
\newblock {\em Advances in Combinatorics}, {J}ul 15 2024.

\bibitem{jiang2024tree}
T.~Jiang, S.~Longbrake, S.~Spiro, and L.~Yepremyan.
\newblock {T}ree {P}osets: {S}upersaturation, {E}numeration, and {R}andomness.
\newblock {\em arXiv:2406.11999}, 2024.

\bibitem{mckinley2023random}
G.~McKinley and S.~Spiro.
\newblock The {R}andom {T}ur{\'a}n {P}roblem for {T}heta {G}raphs.
\newblock {\em arXiv:2305.16550}, 2023.

\bibitem{milicevic2017sets}
L.~Mili\'{c}evi\'{c}.
\newblock Sets in almost general position.
\newblock {\em Combin. Probab. Comput.}, 26(5):720--745, 2017.

\bibitem{morris2016number}
R.~Morris and D.~Saxton.
\newblock The number of {$C_{2\ell}$}-free graphs.
\newblock {\em Adv. Math.}, 298:534--580, 2016.

\bibitem{nenadov2024number}
R.~Nenadov.
\newblock The number of arcs in {$\mathbb{F}_q^2$} of a given cardinality.
\newblock {\em arXiv:2410.21818}, 2024.

\bibitem{nie2024random}
J.~Nie and S.~Spiro.
\newblock Random {T}ur{\'a}n {P}roblems for {$K_{s,t}$} {E}xpansions.
\newblock {\em arXiv:2412.09367}, 2024.

\bibitem{rochenewton2022arcs}
O.~Roche-Newton and A.~Warren.
\newblock Arcs in {$\Bbb{F}_q^2$}.
\newblock {\em European J. Combin.}, 103:Paper No. 103512, 15, 2022.

\bibitem{roth1951problem}
K.~F. Roth.
\newblock On a problem of {H}eilbronn.
\newblock {\em J. London Math. Soc.}, 26:198--204, 1951.

\bibitem{saxton2015hypergraph}
D.~Saxton and A.~Thomason.
\newblock Hypergraph containers.
\newblock {\em Invent. Math.}, 201(3):925--992, 2015.

\bibitem{suk2023higher}
A.~Suk and J.~Zeng.
\newblock On higher dimensional point sets in general position.
\newblock In {\em 39th {I}nternational {S}ymposium on {C}omputational
  {G}eometry}, volume 258 of {\em LIPIcs. Leibniz Int. Proc. Inform.}, pages
  Art. No. 59, 13. Schloss Dagstuhl. Leibniz-Zent. Inform., Wadern, 2023.

\end{thebibliography}
\end{document}